\documentclass[11pt,reqno]{amsproc}
\usepackage{fullpage}
\usepackage{hyperref}
\usepackage{natbib}
\usepackage{cite}
\usepackage{amsmath}
\usepackage{caption}
\usepackage{wrapfig}
\usepackage{lettrine}
\usepackage{morefloats}
\usepackage{bm}        
\usepackage{amssymb}   
\usepackage{amsmath}   
\usepackage{algorithm}
\usepackage{algorithmic}
\usepackage[T1]{fontenc}
\numberwithin{equation}{section}
\usepackage{MnSymbol}
\usepackage{stmaryrd}
\usepackage{graphicx}
\usepackage{subfigure}
\usepackage{enumerate}

\title{Numerical time integration of lumped parameter systems 
governed by implicit constitutive relations}

\author{S. Karimi\\
Email: skarimi@fullranksoftware.com}

\begin{document}

\begin{abstract}

	Time-integration for lumped parameter systems obeying implicit Bingham-Kelvin constitutive models is studied. The governing system of equations describing the lumped parameter system is a non-linear differential-algebraic equation and needs to be solved numerically. The response of this system is non-smooth and the kinematic variables can not be written in terms of the dynamic variables, explicitly. To gain insight into numerical time-integration of this system, a new time-integration scheme based on the trapezoidal method is derived. This method relies on two independent parameters to adjust for damping and is stable. Numerical examples showcase the performance of the proposed time-integration method and compare it to a benchmark algorithm. Under this scheme, implicit-explicit integration of the governing equations is possible. Using this new method, the limitations of the trapezoidal time-integration methods when applied to a non-smooth differential-algebraic equation are highlighted.        
\end{abstract}

\keywords{time-integration methods, lumped systems, 
	differential/algebraic equation, Bingham model, 
	elasto-viscoplastic, non-smooth dynamics}

\maketitle

\section{Introduction}
\label{Sec:Intro}
Constitutive relations play a major role in modeling response of mechanical systems. 
These constitutive relations often can be written as a relation between kinematic 
and dynamic variables. This relation can be explicit, meaning that
the value of a dynamic variable (such as force or stress) is given in terms of kinematic 
variables (such as displacement or strain). In this case, the 
governing differential equations of the mechanical system will be 
in the form of Ordinary or Partial Differential Equations (ODEs or PDEs). 
However, another possible case is the relation between kinematic and dynamic 
variables to be implicit, hence, the value of dynamic variable might not be a function
(one-to-one) relation in terms of kinematic variables \citep{2003_Rajagopal_AoM}.
For example, one might have kinematic variables in terms of dynamic variables. An 
example of these material models is the Bingham model \citep{1922_Bingham,1966_Peryzna_AAM}. 
This case, renders the governing equations of the mechanical system as Differential-Algebraic 
Equations (DAEs) or Partial Differential-Algebraic Equations (PDAEs) \citep{2010_Darbha_JFI}.
Obviously, these equations may not be solved exactly and one needs to resort to numerical 
approximations.

Historically, numerical methods have been developed for ordinary (partial) 
differential equations. Classical examples of time-integration algorithms
are designed to perform well when applied to ODEs. The efforts by 
researchers in this case has been tremendously successful and forms the 
basis for most simulations in scientific computing. However, numerical
integration of DAEs is still in a state of transition and development. 
The reason why numerical integration of DAEs can be more
challenging than ODEs is that the presence of algebraic variables or 
relations can invoke instabilities \citep{1995_Petzold_DAE,1991_Hairer_book}.
In fact, a numerical time-integrator can be successful in approximating
the solution to ODEs, but to fail in estimating the solution to a DAE.  
Considering the importance of DAEs to mathematical modeling of mechanical
systems such as multi-body dynamics \citep{1997_Schiehlen_MSD,2010_BrulsCardona_JCND}, 
domain decomposition for PDEs \citep{2013_Bursi_CM,2014_Karimi_JCP,2015_Karimi_CMAME}, 
as well as implicit constitutive relations in mechanical systems \citep{2010_Rajagopal_MRC,2010_Darbha_JFI,2012_Prazak_AM} 
and many others, many research endeavors have been centered around development of reliable 
numerical integration of DAEs for various applications. It is worthy to note
that the algebraic constraints or relations in a DAE can be of different types.
The algebraic constraint can be smooth (differentiable) in which case 
classical ODE solvers need minor adjustments to give accurate results. However,
as the case in this paper, the algebraic constraint (the implicit constitutive
relation) is continuous but non-smooth (not differentiable). Therefore, 
numerical time-integration needs to be designed with more care.

\begin{figure}
	\centering
	\includegraphics[clip, scale=0.25]{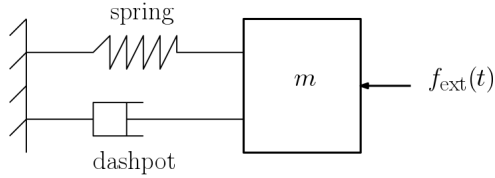}
	\caption{\textsf{SDOF model problem}:~
		In this figure, the SDOF system for the model problem is 
		shown.
	\label{Fig:sdof}}
\end{figure}

In this paper, time-integration of a lumped system (see figure \ref{Fig:sdof}) 
governed by the implicit Bingham-Kelvin constitutive relation is studied. The 
reason why we have chosen this rather simple system is to examine
the basic performance of the derived time-integration methods while avoiding
other numerical complexities. The system considered in this paper, has
a non-smooth response to excitations. Based on findings in \citep{2010_Darbha_JFI}, 
and generalized trapezoidal integration \citep{1990_Wood,2012_Hughes}, 
we derive an integration method that includes two independent parameters. 
A complete description of numerical integration under the proposed methods for 
cases of linear and nonlinear Bingham models is given. Finally,
the performance of the derived methods is compared to the results
from algorithm given in \citep{2010_Darbha_JFI} using several
numerical examples.

\subsection{Outline of the paper}
The rest of the paper is organized as follows. In section \ref{Sec:Method}, 
the governing equations and the new numerical scheme is detailed 
and is followed by numerical examples in section \ref{Sec:Numerical}. 
In section \ref{Sec:Conclusion} the conclusions are drawn and suggestions
for future research are given.

\section{Method for numerical Integration}
\label{Sec:Method}
\subsection{Model problem and governing equations}
\begin{figure}
	\centering
	\includegraphics[clip, scale=0.4]{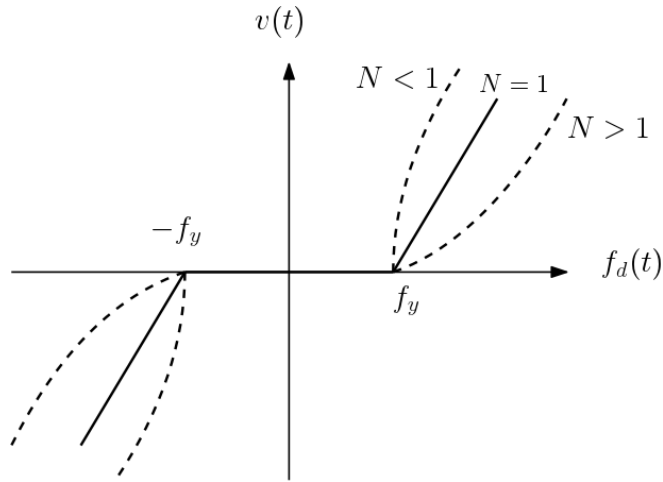}
	\caption{\textsf{Constitutive models of Bingham and Norton}:~
		In this figure, displacement vs. force in Bingham and '
		Norton constitutive relations is shown. For $N = 1$ the 
		Norton model reduces to Bingham model.
	\label{Fig:BinghamNorton}}
\end{figure}

Consider the SDOF system shown in figure \ref{Fig:sdof}, where 
the block of mass $m$ is connected to a linear elastic spring
of stiffness $k$ and a dashpot with Bingham constitutive relation. The 
governing equations for the motion of the block can be written as:
\begin{subequations}
\label{Eqn:SDOF_Equation}
	\begin{align}
		&\dot{v}(t) = \frac{1}{m}\left( f_{\mathrm{ext}}(t)
		- f_{s}(t) - f_{d}(t)\right) \\
		&f_{s}(t) = k u(t) \\
		&v(t) = \varphi(f_{d}(t)) \\
		&\dot{u}(t) = v(t).
	\end{align}
\end{subequations}
where $u(t)$ and $v(t)$ are the displacement and velocity of the
block. The force of dashpot and spring are shown by $f_{d}(t)$ and 
$f_{s}(t)$, respectively. The function $\varphi$
is defined as follows:
\begin{align}
	\label{Eqn:dashpotFun}
	\varphi (f) := \left\{ \begin{array}{l l}
		0 & \big| f \big| \leq f_{y} \\
		\gamma \left( \big| f \big| - f_{y}\right)^{N} 
		\mathrm{Sign}[f] & \big| f \big| \geq f_{y}
		\end{array}
	\right.,
\end{align}
where $\gamma$, $N$ and $f_{y}$ are constitutive parameters. 
Equations \eqref{Eqn:SDOF_Equation} and \eqref{Eqn:dashpotFun}
form a DAE system of equations. In this problem, the velocity is 
given as a function of the force of dashpot, hence, it is referred 
to as an implicit constitutive relation. Figure \ref{Fig:BinghamNorton} 
shows the relation between the dashpot force and the velocity for both 
$N = 1$ and $N \neq 1$. From figure \ref{Fig:BinghamNorton} it can be 
observed that the function $\varphi(f_{d})$ is a nondecreasing function.
Herein, inspired by \citep{2012_Prazak_AM}, we show that the solution to 
the equation \eqref{Eqn:SDOF_Equation} is unique.

\emph{Proposition}. Consider initial conditions $u(t=0)=u_{0}$ and 
$\dot{u}(t=0)=v_{0}$ and let $\varphi$ be a nondecreasing function, 
then the solution to equation \eqref{Eqn:SDOF_Equation} is unique.

\emph{Proof}. Suppose $u$ and $\tilde{u}$ are both solutions to equation
\eqref{Eqn:SDOF_Equation}, with the initial conditions $u(t=0)=u_{0}$ and 
$\dot{u}(t=0)=v_{0}$. Consider the \emph{arbitrary} time-interval $[0,\delta]$. 
Since the initial conditions for both solutions $u$ and $\tilde{u}$ are 
the same, then $u - \tilde{u} = 0$ for $t = 0$. From equation \eqref{Eqn:SDOF_Equation}
we can derive:
\begin{align}
	\label{Eqn:Proof1}
	\ddot{u} - \ddot{\tilde{u}} + \frac{1}{m}(f_{s} - \tilde{f}_{s})
	+ \frac{1}{m}(f_{d} - \tilde{f}_{d}) = 0,
\end{align}
where $f_{s} = k u$ and $\tilde{f}_{s} = k \tilde{u}$. Respectively,
$f_{d}$ and $\tilde{f}_{d}$ correspond to the dashpot force for 
values of $\dot{u}$ and $\dot{\tilde{u}}$. Multiplying equation \eqref{Eqn:Proof1}
by $2(\dot{u} - \dot{\tilde{u}})$ gives
\begin{align}
	\label{Eqn:Proof2}
	2(\dot{u} - \dot{\tilde{u}})(\ddot{u} - \ddot{\tilde{u}}) + 
	\frac{2}{m}(\dot{u} - \dot{\tilde{u}})(f_{s} - \tilde{f}_{s})
	+ \frac{2}{m}(\dot{u} - \dot{\tilde{u}})(f_{d} - \tilde{f}_{d}) = 0.
\end{align}
Note that the nondecreasing assumption on $\varphi$ of equation \eqref{Eqn:dashpotFun}
gives
\begin{align}
	\label{Eqn:Proof3}
	(\dot{u} - \dot{\tilde{u}})(f_{d} - \tilde{f}_{d}) \geq 0.
\end{align}
Denoting $w = u - \tilde{u}$ and using relations \eqref{Eqn:Proof2} and 
\eqref{Eqn:Proof3} result in 
\begin{align}
	\label{Eqn:Proof4}
	\frac{\mathrm{d} \dot{w}^{2}}{\mathrm{d}t} \leq -\frac{k}{2} \frac{\mathrm{d} w^{2}}{\mathrm{d}t}.
\end{align}
Note that based on definition of $w$ we have $w(t=0) = 0$ and $\dot{w}(t=0)=0$.
Integrating \eqref{Eqn:Proof4} over time interval $[0,\delta]$ gives
\begin{align}
	\label{Eqn:Proof5}
	\left(\dot{w}(t = \delta)\right)^{2} + \frac{k}{2} \left( w(t=\delta)\right)^2 \leq 0.
\end{align}
The only way relation \eqref{Eqn:Proof5} can be satisfied is to have
\begin{align}
	\label{Eqn:Proof6}
	\dot{w}(t = \delta) = 0 \quad \text{and} \quad w(t=\delta) = 0.
\end{align}
This result shows that the value of the solution to \eqref{Eqn:SDOF_Equation}
for a give initial condition is unique.$\Box$

Note that the proposition given above applies to both linear and nonlinear Bingham 
models as in both cases $\varphi(f_{d})$ is a nondecreasing function. 

Here, we will differentiate equation \eqref{Eqn:SDOF_Equation}b 
with respect to time in order to reduce the index of this DAE system.
This reduction in DAE index results in a more stable numerical integration.
Doing so, results in the following system of equations:
\begin{subequations}
	\label{Eqn:SDOF_Equation2}
	\begin{align}
		&\dot{v}(t) = \frac{1}{m}\left( f_{\mathrm{ext}}(t)
		- f_{s}(t) - f_{d}(t)\right) \\
		&\dot{f}_{s}(t) = k v(t) \\
		&v(t) = \varphi(f_{d}(t)).
	\end{align}
\end{subequations}
Arguments on existence and uniqueness of a solution to this system 
can be found in \citep{2010_Darbha_JFI,2012_Prazak_AM}, as well as 
using the uniqueness agument given above. In the following section, 
the efforts in reference \citep{2010_Darbha_JFI} 
will be extended and applied to solve the system given in 
\eqref{Eqn:SDOF_Equation2}.

\subsection{Numerical time integration}
Suppose that the time interval of interest is $(0,T]$ and is 
divided into sub-intervals of equal size (time-step) $\Delta t$. 
The value of variable $\Box$ at time level $n\Delta t$ will be 
shown as $\Box^{(n)}$. We will begin by integrating the linear 
momentum balance equation:
\begin{align}
	\label{Eqn:method1}
	\int_{t^{(n)}}^{t^{(n+1)}} \dot{v}(t) \mathrm{d}t = 
	\frac{1}{m} \int_{t^{(n)}}^{t^{(n+1)}} \left( f_{\mathrm{ext}}(t)
	 - f_{s}(t) - f_{d}(t)\right) \mathrm{d}t.
\end{align}
Using trapezoidal integration, equation \eqref{Eqn:method1} results
in:
\begin{align}
	\label{Eqn:method2}
	v^{(n+1)} - v^{(n)} = \frac{\Delta t}{m} \left( (1 - \alpha)
	\left( f_{\mathrm{ext}}^{(n)} - f_{s}^{(n)} - f_{d}^{(n)}\right) 
	+ \alpha \left( f_{\mathrm{ext}}^{(n+1)} - f_{s}^{(n+1)} - 
	f_{d}^{(n+1)}\right)\right),
\end{align}
where $\alpha$ is the integration parameter ($0\leq\alpha\leq1$). The 
equation \eqref{Eqn:SDOF_Equation2}b can be integrated in a similar
fashion to give:
\begin{align}
	\label{Eqn:method3}
	f_{s}^{(n+1)} - f_{s}^{(n)} = k\Delta t \left( (1-\beta)v^{(n)}
	+ \beta v^{(n+1)}\right),
\end{align}
where $\beta$ is an integration parameter in $[0,1]$. Equations 
\eqref{Eqn:method2} and \eqref{Eqn:method3} can be combined to give:
\begin{align}
	\label{Eqn:method4}
	\left( 1 + \frac{\alpha \beta \Delta t^{2} k}{m}\right) v^{(n+1)}
	= \frac{\alpha \Delta t}{m}\left( \hat{f}^{(n+1)} - f_{d}^{(n+1)}\right),
\end{align}
where $\hat{f}^{(n+1)}$ is defined as:
\begin{align}
	\label{Eqn:method5}
	\hat{f}^{(n+1)} = f_{\mathrm{ext}}^{(n+1)} + 
	\left( \frac{1}{\alpha} - 1\right)f_{\mathrm{ext}}^{(n)} - 
	\frac{1}{\alpha} f_{s}^{(n)} - \left( \frac{1}{\alpha} - 1\right)f_{d}^{(n)}
	+ \left( \frac{m}{\alpha \Delta t} - k \Delta t (1 - \beta)\right)v^{(n)}.
\end{align}
The constitutive relation for dashpot given in equation \eqref{Eqn:dashpotFun}
shows us that $\mathrm{Sign}[v] = \mathrm{Sign}[f_{d}]$ at all time levels.
Also, from equation \eqref{Eqn:method4} we can derive that $\mathrm{Sign}[v^{(n+1)}]
= \mathrm{Sign}[\hat{f}^{(n+1)}]$, because
\begin{align}
	\label{Eqn:method6}
	\underbrace{\left(\left( 1 + \frac{\alpha \beta \Delta t^{2} k}{m}\right) 
	\big|v^{(n+1)}\big| + \frac{\alpha \Delta t}{m} \big| f_{d}^{(n+1)}\big| 
	\right)}_{\geq 0}\mathrm{Sign}[v^{(n+1)}]
	= \underbrace{\frac{\alpha \Delta t}{m} \big| \hat{f}^{(n+1)}}_{\geq 0} 
	\big| \mathrm{Sign}[\hat{f}^{(n+1)}].
\end{align}
Equation \eqref{Eqn:method6} indicates that the value of $\hat{f}^{(n+1)}$
can be used as a predictor to determine the relation between $v^{(n+1)}$
and $f_{d}^{(n+1)}$. This equation implies that if $\big| f_{d}^{(n+1)} \big| \geq f_{y}$,
then we also have $\big| \hat{f}^{(n+1)} \big| \geq \big| f_{d}^{(n+1)}\big|$.
Therefore, if $\big| \hat{f}^{(n+1)} \big| \leq f_{y}$
then $v^{(n+1)} = 0$ and from equation \eqref{Eqn:method4} we get 
$f_{d}^{(n+1)} = \hat{f}^{(n+1)}$. Otherwise, the force of dashpot needs
to be solved for using the following equation:
\begin{align}
	\label{Eqn:method7}
	\left( 1 + \frac{\alpha \beta \Delta t^{2} k}{m}\right)\gamma
	\left( \big| f_{d}^{(n+1)} \big| - f_{y}\right)^{N} \mathrm{Sign}[\hat{f}^{(n+1)}]
	- \frac{\alpha \Delta t}{m}\left( \hat{f}^{(n+1)} - f_{d}^{(n+1)}\right)=0.
\end{align}
Equation \eqref{Eqn:method7} is obtained from equations \eqref{Eqn:method4}
and \eqref{Eqn:dashpotFun}. Obviously, this equation is in general a 
nonlinear equation. In the case of $N = 1$, equation
\eqref{Eqn:method7} can be simplified as:
\begin{align}
	\label{Eqn:method8}
	f_{d}^{(n+1)} = \frac{1}{1 + \frac{\alpha}{\gamma m \left( \frac{1}{\Delta t}
	+ \frac{\alpha \beta \Delta t k}{m}\right)}} \left( \frac{\alpha}{\gamma m 
	\left( \frac{1}{\Delta t} + \frac{\alpha \beta \Delta t k}{m}\right)} 
	\hat{f}^{(n+1)} + \mathrm{Sign}[\hat{f}^{(n+1)}] f_{y}\right).
\end{align}
Hence a nonlinear solver is not required. Once the value of the $f_{d}^{(n+1)}$ 
is found, the value of $v^{(n+1)}$ can be found to be 
\begin{align}
	\label{Eqn:method9}
	v^{(n+1)} = \gamma \left( \big| f_{d}^{(n+1)} \big| - f_{y}\right)^{N} 
	\mathrm{Sign}[f_{d}^{(n+1)}].
\end{align}
The value of the spring force can be updated using equation \eqref{Eqn:method3} 
and the displacement of the mass can be found from the following relation:
\begin{align}
	\label{Eqn:method10}
	u^{(n+1)} = \frac{1}{k} f_{s}^{(n+1)}.
\end{align}
To complete the description of time-stepping algorithm, we need to clarify
two points:
\begin{enumerate}[(a)]
\item \textsf{Initialization}:~In this paper, we will assume all the 
	initial displacements to be purely elastic. In other words, given an
	initial displacement $u_{0}$, the initial spring force will be taken
	to be $f_{s,0} = k u_{0}$ and initial dashpot force $f_{d,0}$ should be 
	found from
	\begin{align}
		\label{Eqn:method11}
		\big| v_{0} \big| = \gamma \left( \big| f_{d,0}\big| - f_{y}\right)^{N}
	\end{align}
	if initial velocity $v_{0}$ is nonzero. Note that $\mathrm{Sign}[f_{d,0}] = 
	\mathrm{Sign}[v_{0}]$. Consequently, equation \eqref{Eqn:method11} for the 
	case of $N = 1$ gives:
	\begin{align}
		\label{Eqn:method12}
		f_{d,0} = \left( f_{y} + \frac{\big| v_{0} \big|}{\gamma} \right) \mathrm{Sign}[v_{0}].
	\end{align}

\item \textsf{Initial guess for nonlinear solver}:~Equation \eqref{Eqn:method7}
	is a nonlinear equation with respect to $f_{d}^{(n+1)}$. Hence, as with 
	any nonlinear solver, an initial guess is needed. Based on our numerical 
	experiments, the initial guess of $f_{y} \mathrm{Sign}[\hat{f}^{(n+1)}]$
	seems to be a reasonable choice.
	
\end{enumerate}
This concludes the description of the new time-stepping algorithm. A pseudocode
for this method is given in algorithm \ref{Alg:NumMethod}. In the 
following section, numerical examples will be presented to demonstrate
the effect of various time-steps and time-integration parameters $\alpha$
and $\beta$.
\begin{algorithm}
	\caption{\textsf{Numerical time-integration algorithm}:~The pseudocode 
	for time-stepping algorithm for elstic-viscoplastic SDOF system.}
	\label{Alg:NumMethod}
	\begin{algorithmic}[1]
		\STATE Initialize $f_{s,0}$, $f_{d,0}$ for given $u_{0}$ and $v_{0}$. 
		\STATE Set values for $\Delta t$, $\alpha$ and $\beta$.
		\STATE Set $t = 0$.
		\WHILE{$t \leq T$} 
			\STATE Set $t = t + \Delta t$.
			\STATE Calculate the predictor $\hat{f}$ from equation \eqref{Eqn:method5}.
			\IF{$\big| \hat{f} \big| \leq f_{y}$}
				\STATE Set the value of velocity $v = 0$.
				\STATE Set $f_{d} = \hat{f}$.
			\ELSE
				\STATE Find the value of dashpot force $f_{d}$ from equation \eqref{Eqn:method7}.
				\STATE Find the value of velocity from equation \eqref{Eqn:method9}.
			\ENDIF
			\STATE Find the value of the spring force $f_{s}$ from equation \eqref{Eqn:method3}.
			\STATE Update the value of displacement from equation \eqref{Eqn:method10}.
		\ENDWHILE
	\end{algorithmic}
\end{algorithm}

\section{Numerical Results}
\label{Sec:Numerical}
In this section, we will showcase the performance of the derived method in 
comparison to the method proposed in \citep{2010_Darbha_JFI}. The amount of
damped (dissipated) energy will be denoted by $E_{d}(t)$ and is defined as
\begin{align}
	\label{Eqn:Edissip}
	E_{d}(t) := \int_{0}^{t} v(t) f_{d}(t) \;\mathrm{d}t.
\end{align}
Note that because $v(t)$ and $f_{d}(t)$ have the same sign, the value
of $E_{d}(t)$ will be nonnegative at all time-levels. The error in
kinematic quantities will be denoted by $e_{p}$ for $p=u,v$, and is
defined as
\begin{align}
	\label{Eqn:error}
	e_{p} := \frac{1}{M}\sqrt{\sum_{i=1}^{M}(p^{(i)} - p_{\mathrm{ref}}^{(i)})^{2}}\;,
\end{align}
where $M$ is the total number of time-steps and $p_{\mathrm{ref}}$ is
the reference (benchmark) values. The numerical results are obtained
from implementation in MATLAB \citep{MATLAB}.
\subsection{The case of $N = 1$}
In this problem, consider the SDOF system shown in figure \ref{Fig:sdof}
subject to external force given as:
\begin{align}
	\label{Eqn:fext}
	f_{\mathrm{ext}}(t) = 2 \sin (2 \pi t) e^{(-0.2 t)}.
\end{align}
The initial conditions are $u_{0} = 0$ and $v_{0} = 0$. The physical parameters 
for this problem are given in Table \ref{Tbl:PhysBingN=1}.
To demonstrate the effect of time-integration parameters $\alpha$ and $\beta$,
we use different combinations of values of either 1 or 1/2 for each of them.
Time-integration parameters used in all cases are given in Table \ref{Tbl:BingN=1}.
Note that the benchmark in Table \ref{Tbl:BingN=1} refers to the method
proposed in \citep{2010_Darbha_JFI}, which coincides with $\alpha=1$ and $\beta=1$
under the method given in this paper. All the numerical solutions will be compared
to that of benchmark method.

\begin{figure}
	\centering
	\subfigure[displacement]{\includegraphics[scale=0.33,clip]{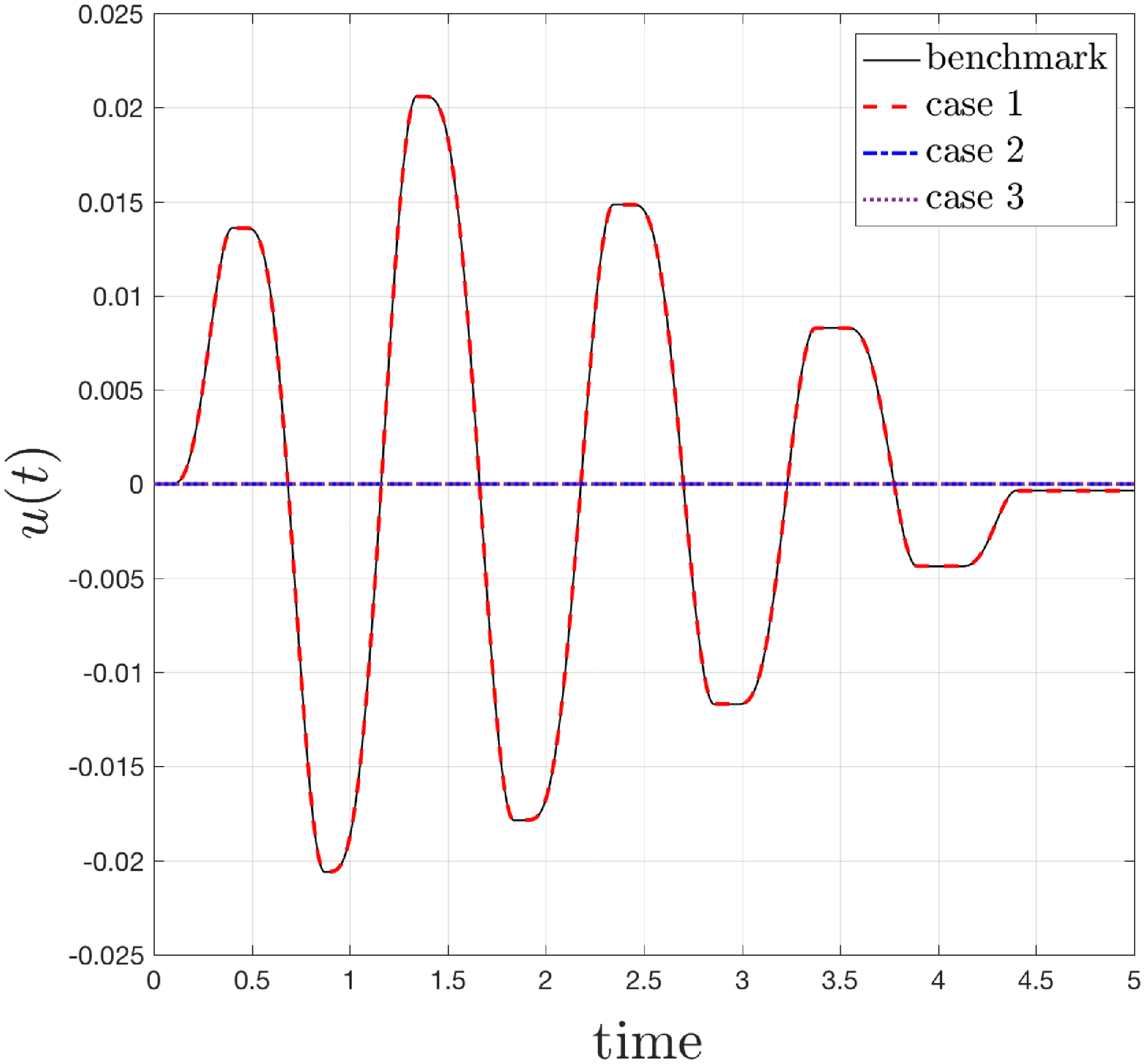}}
	\subfigure[velocity]{\includegraphics[scale=0.33,clip]{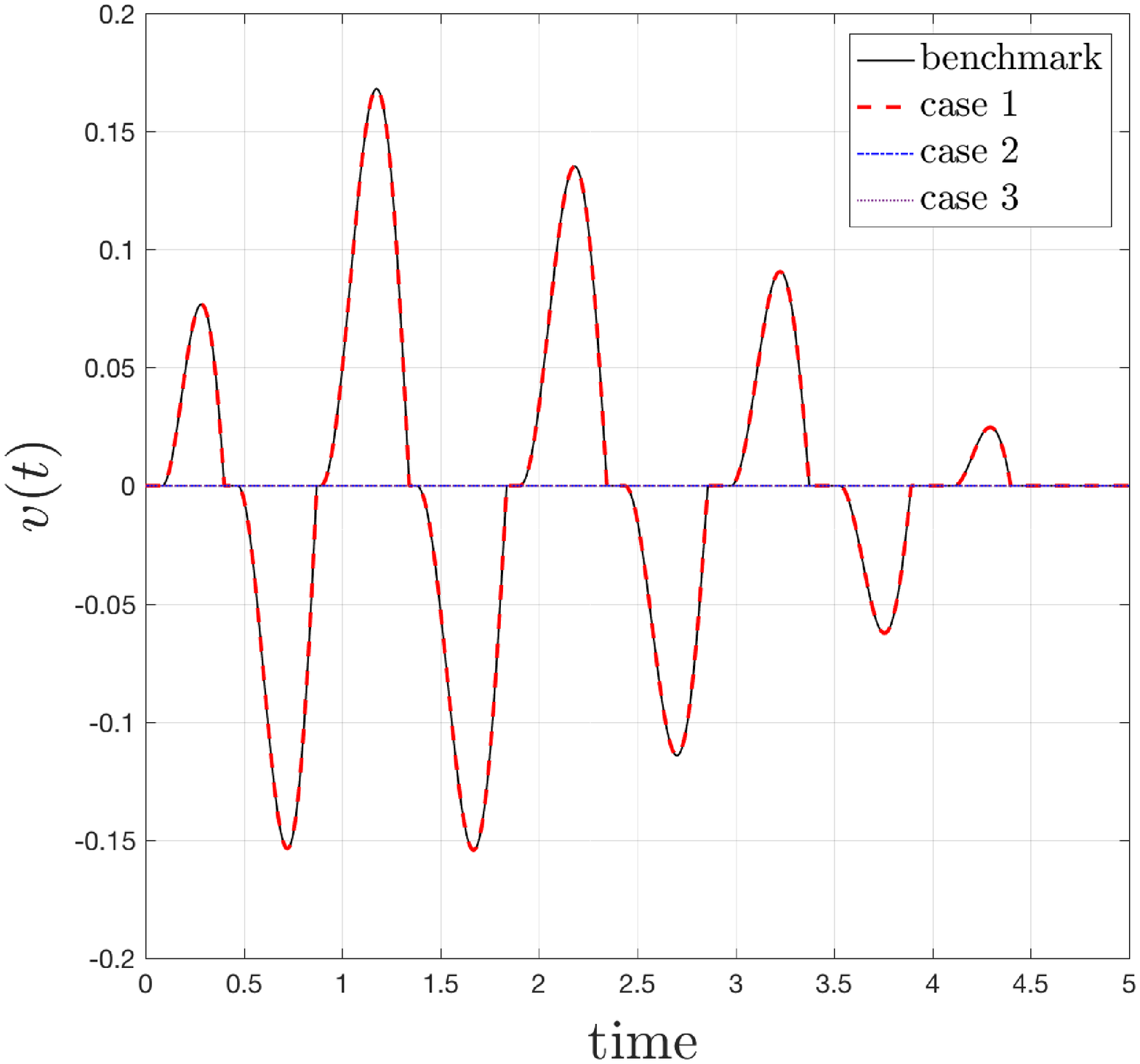}}
	\caption{\textsf{Displacement and velocity for case of $N=1$}:~ In this figure, the 
	displacement and velocity against time are shown. As it can be seen, only the answer 
	from the case $\alpha = 1$ and $\beta = 1/2$ are comparable to the benchmark solution.
	\label{Fig:Bing_DispVel}}
\end{figure}

\begin{figure}
	\centering
	\subfigure[force of spring]{\includegraphics[scale=0.33,clip]{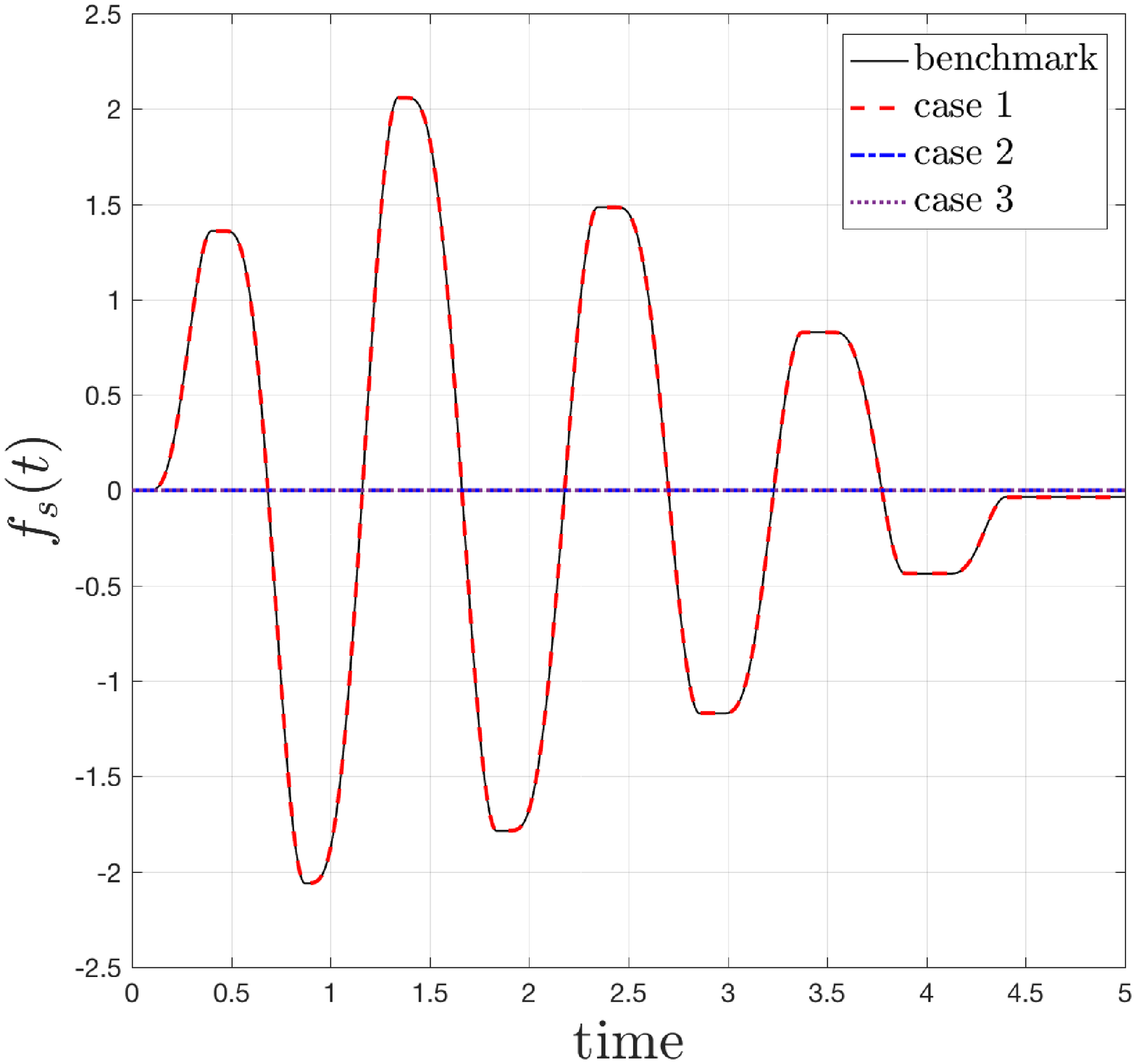}}
	\subfigure[force of dashpot]{\includegraphics[scale=0.33,clip]{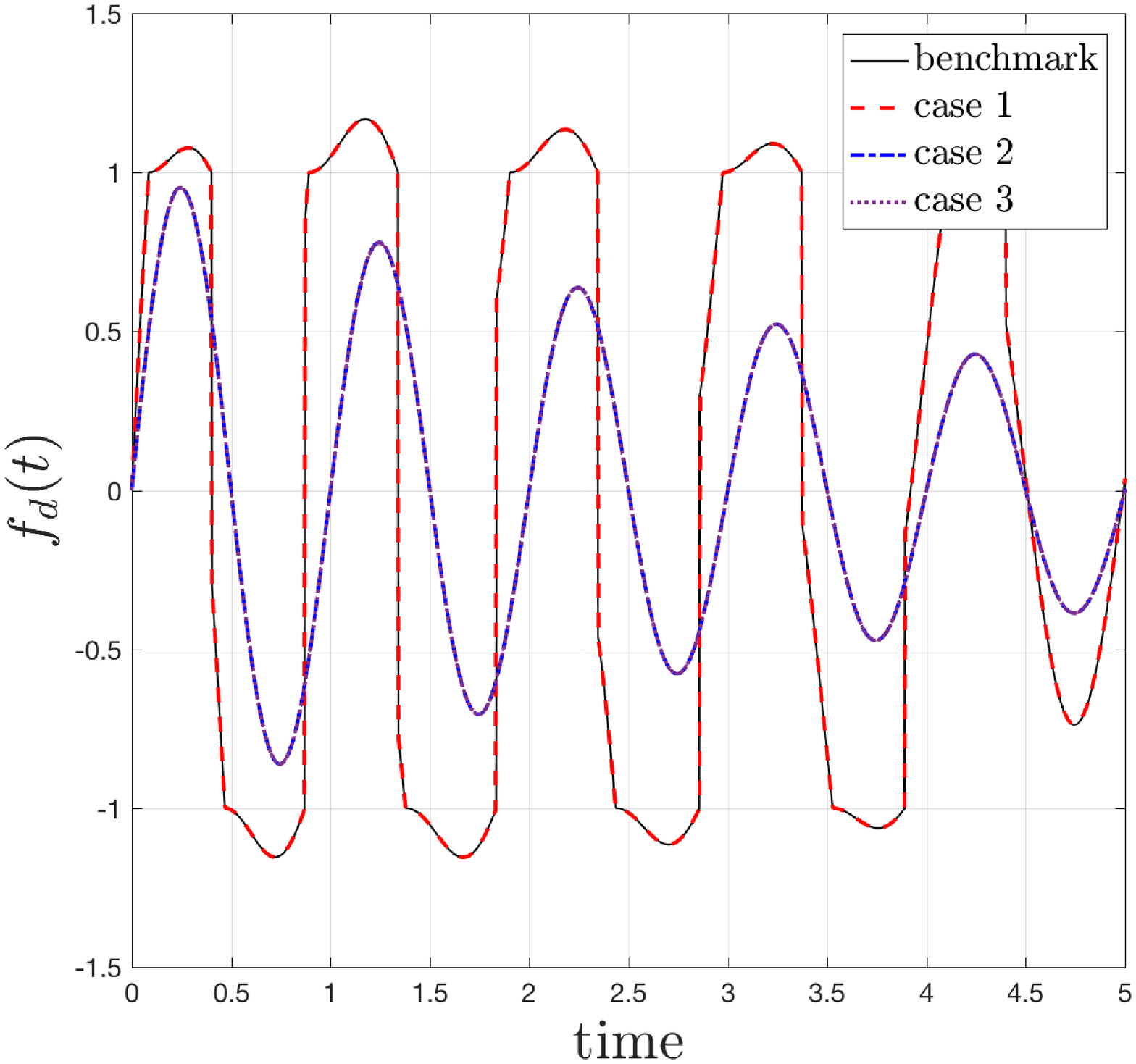}}
	\caption{\textsf{Spring and dashpot force for case of $N=1$}:~The force of the spring
	and dashpot are shown against time. In cases 2 and 3, the damping in the magnitude of the 
	dashpot force results in inaccurate solutions for velocity, displacement and the spring
	force $f_{s}(t)$.
	\label{Fig:Bing_force}}
\end{figure}

\begin{figure}
	\centering
	\subfigure[$\Delta t = 10^{-6}$]{\includegraphics[clip,scale=0.4]{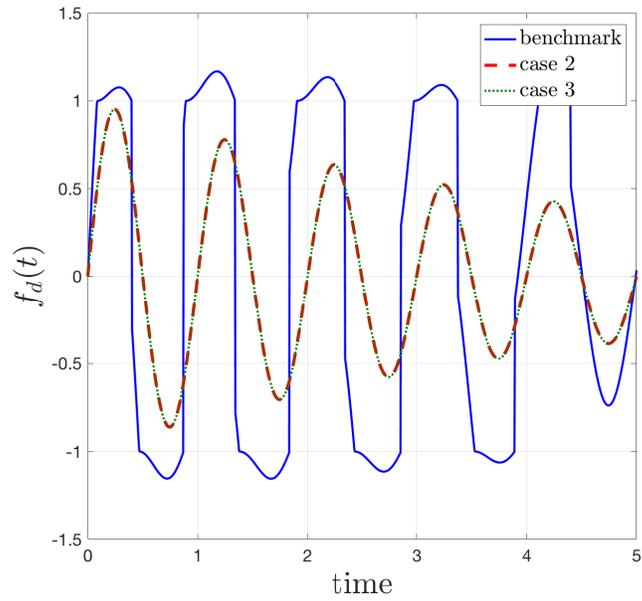}}
	\subfigure[$\Delta t = 10^{-7}$]{\includegraphics[clip,scale=0.4]{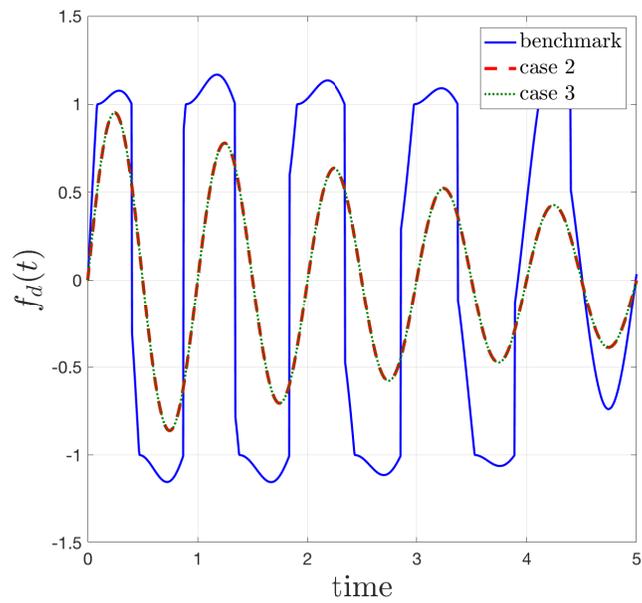}}
	\caption{\textsf{Dashpot force for $N=1$ with smaller time-steps}:~
	In this figure, the force of dashpot compared to the solution from the benchmark
	method is shown. The values of $\alpha$ and $\beta$ are the same as in cases 2 and 
	3 of Table \ref{Tbl:BingN=1}. It can be seen that refining time-step does not improve 
	damping properties of these methods. 
	\label{Fig:Bingham_case23}}
\end{figure}

\begin{figure}
	\centering
	\includegraphics[clip, scale=0.4]{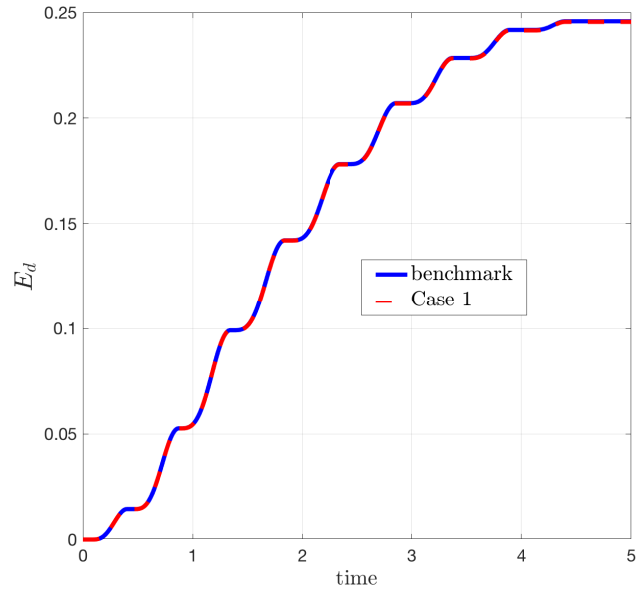}
	\caption{\textsf{Energy dissipation}:~ The amount of energy dissipated over time is shown 
	in this figure. The values from the case 1 and the benchmark closely match. 
	\label{Fig:BingDissip}}
\end{figure}

\begin{figure}
	\centering
	\includegraphics[scale=0.4,clip]{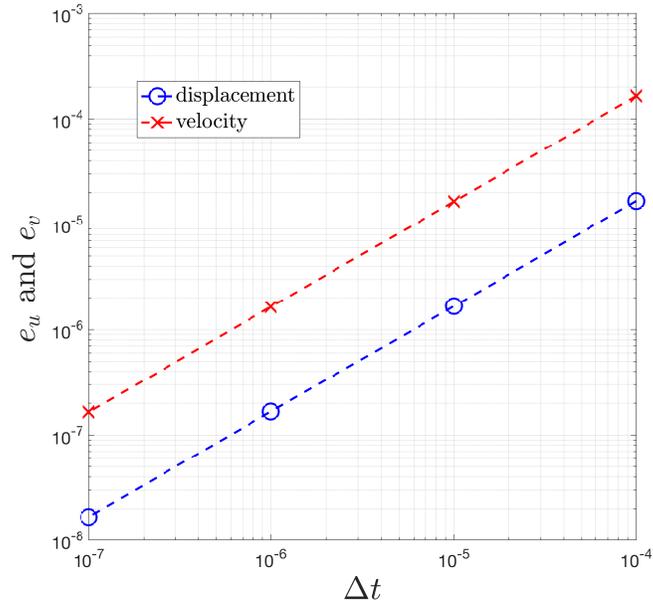}
	\caption{\textsf{Error compared to the benchmark for the case $N=1$}:~
	This figure shows the error in velocity and displacement compared to the 
	numerical solution of the benchmark method proposed in \citep{2010_Darbha_JFI}.
	It can be seen that this error reduces with respect to time-step in 
	$\mathcal{O}(\Delta t)$.
	\label{Fig:Bing_error}}
\end{figure}

The numerical results for displacement and velocity, against time 
are shown in figures \ref{Fig:Bing_DispVel}. The dashpot 
force $f_{d}$ and spring force $f_{s}$ is shown in figure \ref{Fig:Bing_force}. 
It can be seen that the numerical results from the benchmark solution and case 1 
(i.e., $\alpha=1,\;\beta=1$ and $\alpha=1,\;\beta=1/2$) are very different from the 
cases 2 and 3. The numerical solution from case 1 is indeed very close to the solution 
from the benchmark method. Repeating the same simulation with cases 2 and 3 of Table 
\ref{Tbl:BingN=1} with smaller time-steps show that the damping properties in 
values of dashpot force $f_{d}$ do not improve (see figure \ref{Fig:Bingham_case23}). 
In fact the numerical values are very close to the ones with time-step of $10^{-4}$. 
An important physical quantity is the energy damping in equation \eqref{Eqn:Edissip}, 
which shown in figure \ref{Fig:BingDissip} for benchmark case and case 1 of Table \ref{Tbl:BingN=1}. 
The amount of energy dissipated calculated from each of these cases are very similar, 
while cases 3 and 4 suggest zero energy damping. The solution from the benchmark 
method and case 1 converge to one another in order $\mathcal{O}(\Delta t)$, 
shown in figure \ref{Fig:Bing_error}. However, cases 2 and 3 underestimate 
the dashpot force $f_{d}$ resulting in inaccurate values for displacement and 
velocity. In fact, the value of $f_{d}$ in cases 2 and 3 seems to dampen excessively. 

From this numerical experiment, we find that among the choices for parameters $\alpha$
and $\beta$, the values for case 1 and the benchmark method are the only viable ones.
It seems that choosing $\alpha = 1/2$, introduces excessive damping to the value of 
$f_{d}$. Due to the non-smooth nature of the constitutive relation considered in this paper,
this property results in wrong values for other variables. 

\begin{table}
	\centering
	\caption{\textsf{Parameters for case of $N=1$}:~Values for physical parameters
	are given.}
	\label{Tbl:PhysBingN=1}
	\begin{tabular}{ c c c c c c } \hline
	Parameter & $m$ & $k$ & $f_{y}$ & $\gamma$ & $N$ \\ \hline
	Value & 1 & 100 & 1 & 1 & 1 \\ \hline
	\end{tabular}
\end{table}
\begin{table}
	\centering
	\caption{\textsf{Parameters for case of $N=1$}:~The time-integration
	parameters for various cases are given.}
	\label{Tbl:BingN=1}
	\begin{tabular}{ c c c c }\hline
		Case & $\Delta t$ & $\alpha$ & $\beta$ \\ \hline
		benchmark & $10^{-6}$ & 1 &  1 \\
		1 & $10^{-4}$ & 1 & 1/2 \\
		2 & $10^{-4}$ & 1/2 & 1 \\
		3 & $10^{-4}$ & 1/2 & 1/2 \\ \hline
	\end{tabular}
\end{table}

\subsection{The case of $N \neq 1$}
In this case, the value of dashpot force $f_{d}$ is the solution of the nonlinear
equation \eqref{Eqn:method7}. For this case, the physical parameters are given in
Table \ref{Tbl:PhysNort}. Time-integration parameters are given in Table \eqref{Tbl:Nort}.
Here, we take the value of $N = 3$ and the benchmark numerical solution is from the
method in \citep{2010_Darbha_JFI}. The external force is the same as in equation 
\eqref{Eqn:fext}, and the initial displacement and velocity are both zero.

\begin{figure}
	\centering
	\subfigure[displacement]{\includegraphics[scale=0.33,clip]{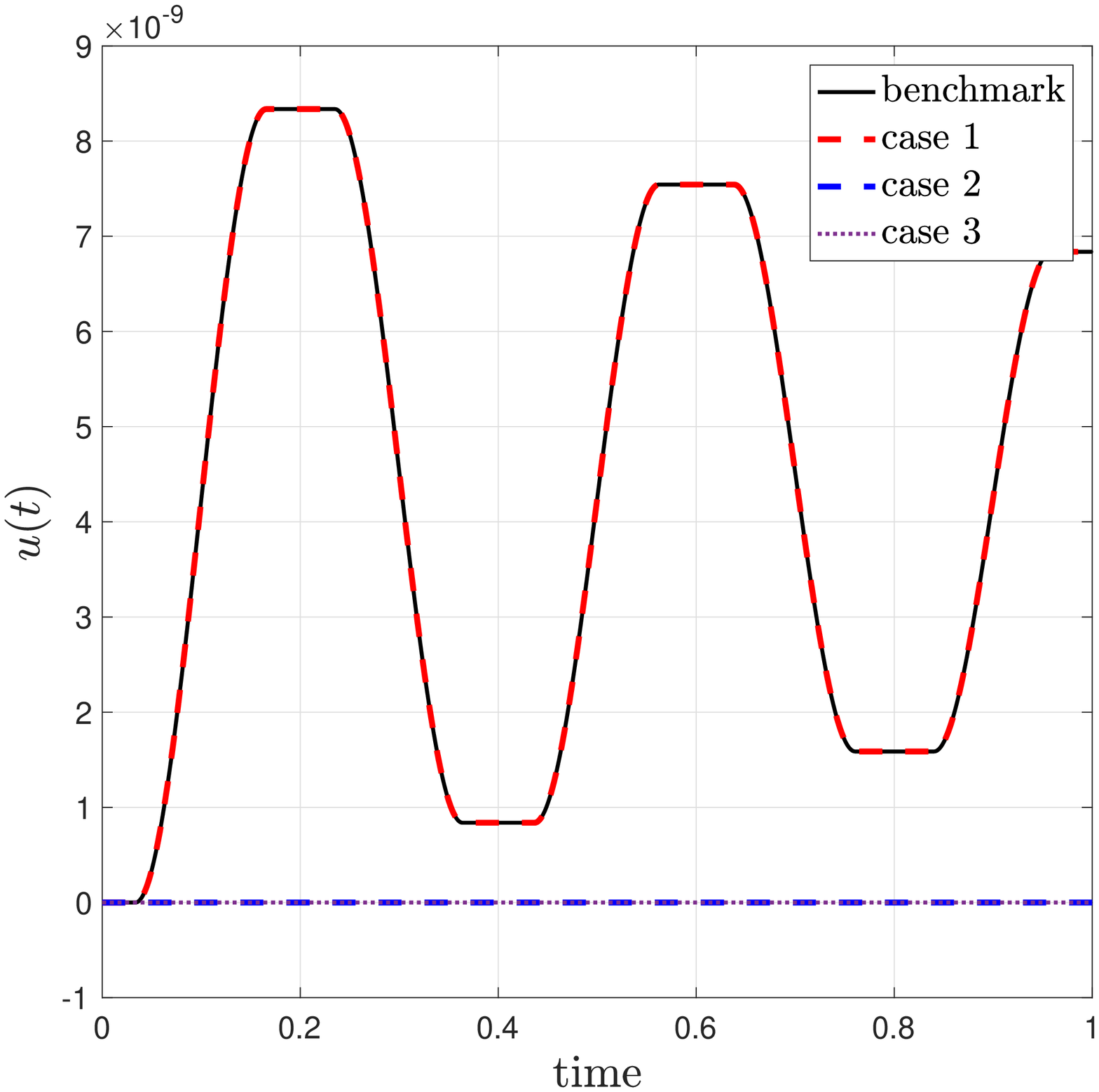}}
	\subfigure[velocity]{\includegraphics[scale=0.33,clip]{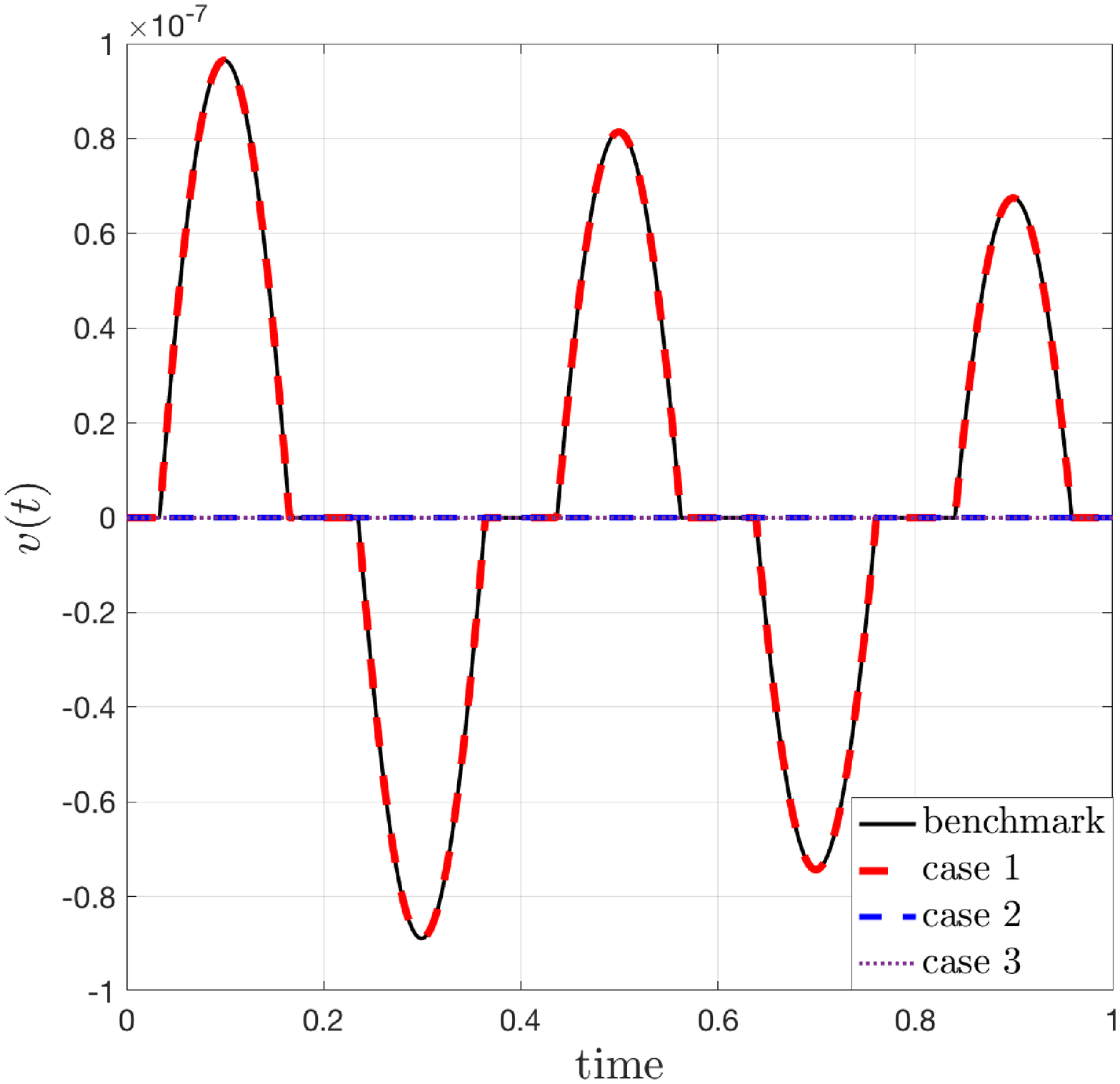}}
	\caption{\textsf{Displacement and velocity for case of $N\neq1$}:~ Comparing the numerical
	results from different cases leads to the conclusion that the only viable choice for 
	parameter $\alpha$ is 1. The only set of results close to the solution from the benchmark
	method is case 1.
	\label{Fig:Nort_DispVel}}
\end{figure}

\begin{figure}
	\centering
	\subfigure[force of spring]{\includegraphics[scale=0.33,clip]{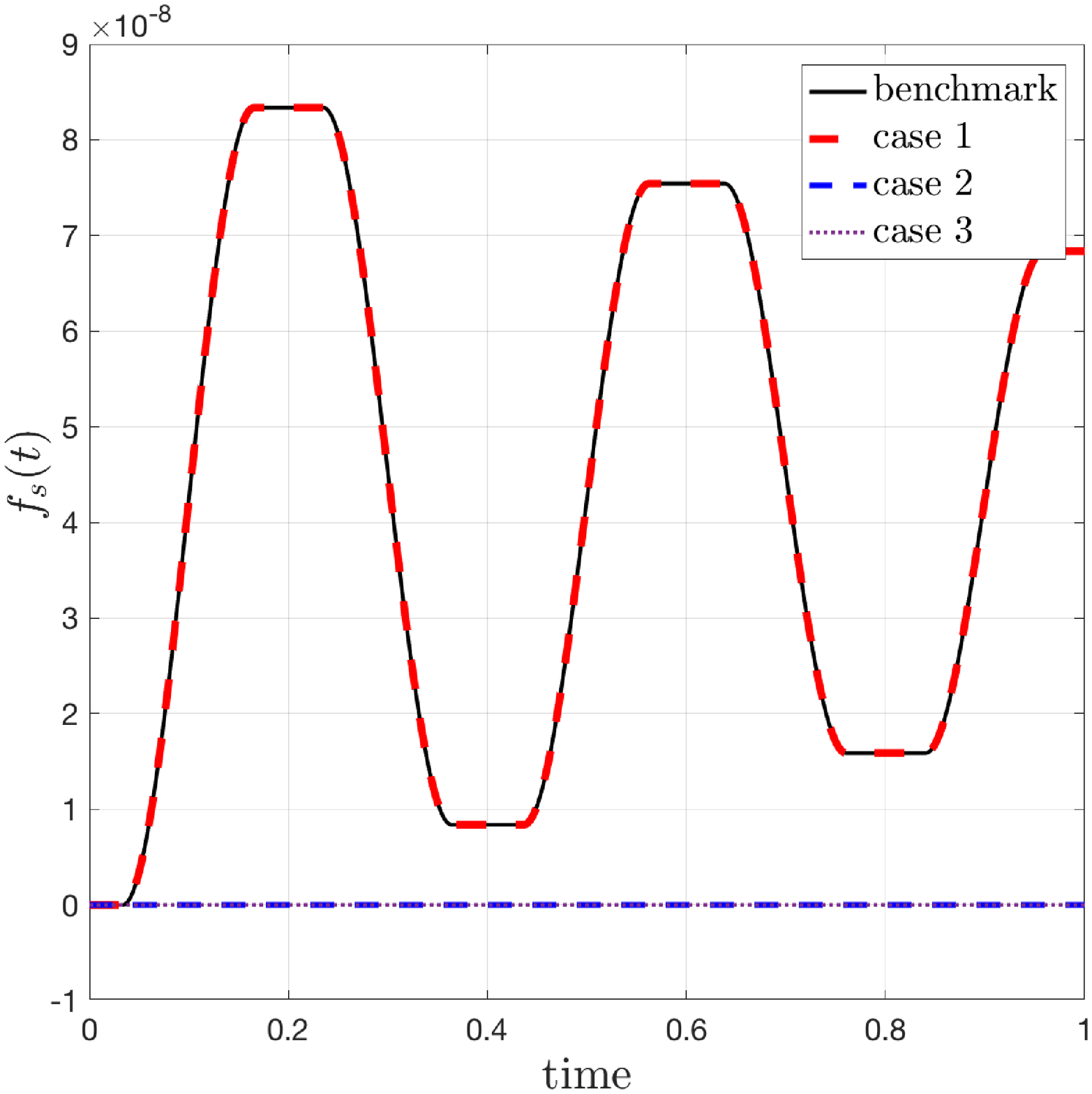}}
	\subfigure[force of dashpot]{\includegraphics[scale=0.33,clip]{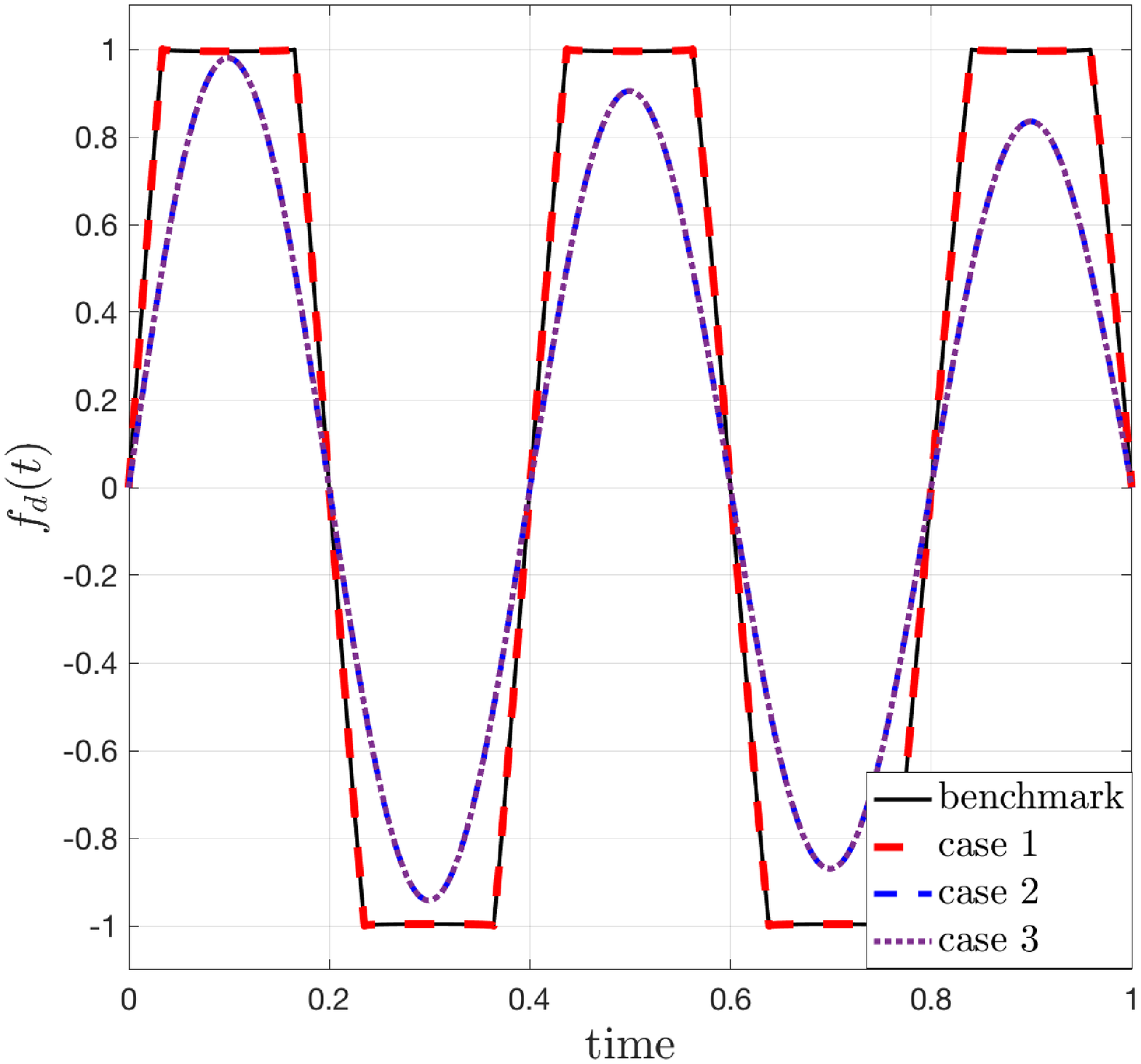}}
	\caption{\textsf{Spring and dashpot force for case of $N\neq1$}:~ The force in spring and 
	dashpot are shown against time. Among all cases, only case 1 gives results that are comparable
	to the benchmark method. Other cases show excessive damping in value of dashpot force.
	\label{Fig:Nort_force}}
\end{figure}

\begin{figure}
	\centering
	\includegraphics[clip, scale=0.4]{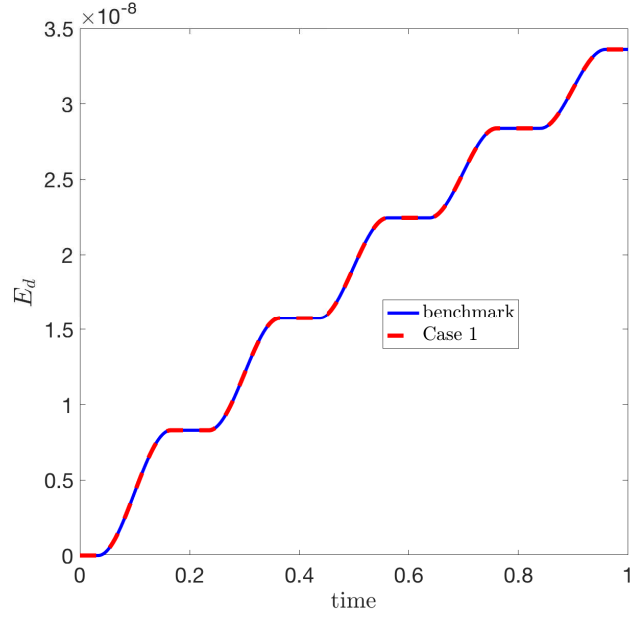}
	\caption{\textsf{Energy dissipation}: The energy dissipation is shown. The 
	value of total amount of energy dissipated by the dashpot is close to the 
	values from the benchmark method. \label{Fig:NortDissip}}
\end{figure}

The numerical results are given in figures \ref{Fig:Nort_DispVel}--\ref{Fig:Nort_force}.
Similar to the case of $N = 1$, the benchmark numerical solution and case 1 are 
close and compatible with one another. However, the parameters selected for cases 2 and 
3 are very different and predict a value of zero throughout the time-interval of interest
for velocity and displacement. In fact, in cases where $\alpha = 1/2$ the smoothing 
results in damping in values of $f_{d}$. It is worth noting that unlike the case of $N=1$,
the displacement seems to converge to a point other than $u_{0}$ as time proceeds.
In the case of $N\neq1$, the amount of energy damped shown in figure \ref{Fig:NortDissip} 
suggests a similar conclusion to the case of $N=1$. In this example, the energy dissipation
predicted by case 1 is very similar to the benchmark solution. Other cases of $\alpha$
and $\beta$ predict zero energy dissipation in this system.

This experiment concludes that the only viable value for $\alpha$ under the proposed
method is 1. The numerical solution to this nonlinear and non-smooth system is sensitive
to smoothing effects of trapezoidal integration in linear momentum balance equation. 
The numerical results seem to be much less sensitive to the choice of $\beta$ between 
1/2 and 1.

\begin{table}
	\centering
	\caption{\textsf{Parameters for case of $N\neq1$}:~Physical
	parameters in the case of $N \neq 1$ are given.}
	\label{Tbl:PhysNort}
	\begin{tabular}{ c c c c c c } \hline
	Parameter & $m$ & $k$ & $f_{y}$ & $\gamma$ & $N$ \\ \hline
	Value & 1 & 10 & 1 & 1 & 3 \\ \hline
	\end{tabular}
\end{table}
\begin{table}
	\centering
	\caption{\textsf{Parameters for case of $N\neq1$}:~Values of the 
	time-integration parameters.}
	\label{Tbl:Nort}
	\begin{tabular}{ c c c c }\hline
		Case & $\Delta t$ & $\alpha$ & $\beta$ \\ \hline
		benchmark & $10^{-7}$ & 1 &  1 \\
		1 & $10^{-7}$ & 1 & 1/2 \\
		2 & $10^{-7}$ & 1/2 & 1 \\
		3 & $10^{-7}$ & 1/2 & 1/2 \\ \hline
	\end{tabular}
\end{table}

\subsection{Implicit-explicit time-integration}
Thus far, the numerical experiments suggest that the only viable value for 
$\alpha$ is 1. However, for $\alpha=1$, the numerical results for various 
values of $\beta$ seem to be very close. In this section, we showcase the 
numerical results for the case of $\alpha=1$ and $\beta = 0$. For $\beta=0$,
the equation \eqref{Eqn:method3} will become an explicit Euler integration
scheme while equation \eqref{Eqn:method2} remains implicit Euler integration.
In this case, the constitutive relation of the elastic spring is integrated
using an explicit scheme and the momentum balance is integrated via an 
implicit relation.

\begin{figure}
	\centering
	\includegraphics[clip, scale=0.4]{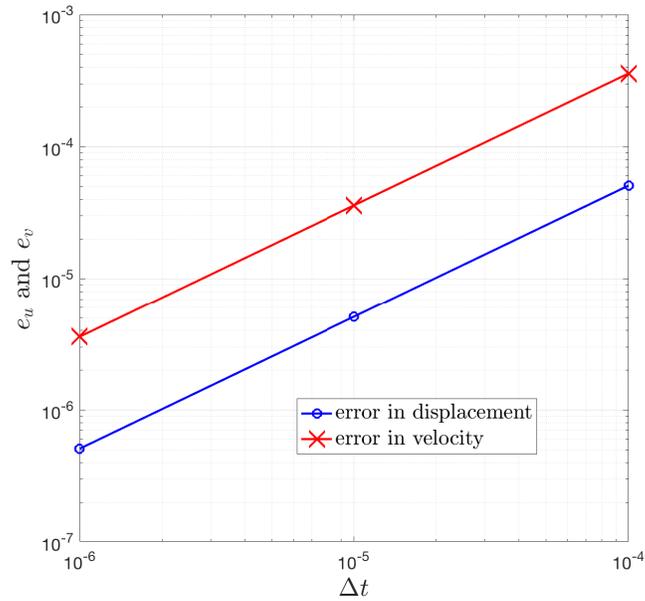}
	\caption{\textsf{Implicit-explicit time-integration}: In this figure the error
	in the implicit-explicit integration ($\alpha=1$ and $\beta=0$) compared to the 
	benchmark numerical method is shown. The numerical values for displacement and 
	velocity converge to the benchmark result in $\mathcal{O}(\Delta t)$. 
	\label{Fig:ImExError}}
\end{figure}

The problem is solved for the external force in equation \eqref{Eqn:fext} using
time-steps $\Delta t = 10^{-4}, 10^{-5}$ and $10^{-6}$. The error (compared to
the benchmark numerical result) is shown in figure \ref{Fig:ImExError}. It can be
seen that the numerical solutions from the two methods converge as $\mathcal{O}(\Delta t)$, 
which implies that the choice of $\alpha=1$ and $\beta=0$ is a viable choice of 
time-integration parameters. This method, therefore, does not suffer from excessive 
numerical damping for values of dashpot force $f_{d}$.

\section{Conclusion}
\label{Sec:Conclusion}
In this paper the governing differential-algebraic equations for a 
SDOF system obeying Bingham constitutive model were reviewed and 
a proof for uniqueness of the solution was presented. Using the 
generalized trapezoidal integration a method for numerical integration
of the given equations was derived with two independent time-integration 
parameters. A pseudocode for this numerical
algorithm along with initialization and initial guess for the 
nonlinear equations was provided. This derivation was followed by 
numerical examples showcasing its performance compared to a
benchmark method. The findings from these numerical experiments 
can be summarized as below:
\begin{enumerate}[(1)]
	\item Because of the non-smooth nature of the physical system 
	considered in this paper, accuracy in estimation of dashpot force
	is of utmost importance. Choice of time-integration parameters 
	can lead to excessive numerical damping in value of dashpot force
	, which in turn gives completely different numerical values for
	other kinematic and dynamic variables.
	\item Under the presented numerical scheme in this paper, the 
	only viable choice of integration parameter $\alpha$ is 1. 
	Other values of $\alpha$ result in additional (unphysical) 
	smoothing in the numerical solution. 
	\item Using the method of numerical time integration in this paper
	the value of integration parameter $\beta$ can be anything in range 
	$[0,1]$. It seems that the values of $\beta$ do not introduce any 
	more significant damping to the estimated value of dashpot force.
	\item Given $\alpha=1$ and $\beta \in [0,1]$ under the numerical 
	scheme here, the numerical solution converges to the solution 
	from the benchmark method in order $\mathcal{O}(\Delta t)$.
	\item Under the proposed method in this paper, one can use the 
	$\alpha=1$ and $\beta=0$ values for integration parameters. This
	means that one can have a convergent implicit-explicit time-integration 
	of the non-smooth system considered here.
\end{enumerate}
Future developments in this area can be extension of such methods to 
other implicit and non-smooth material models Maxwell type
\citep{ChristensenBook}, extension to 2 and 3-dimensional problems, 
and extension of Newmark-type integration methods \citep{1959_Newmark,1993_Hulbert_JAM} for this class of problems.

\bibliographystyle{plainnat}
\bibliography{References}




\end{document}